\documentclass[10pt]{elsarticle}

\usepackage{amsfonts}
\usepackage{xcolor}
\usepackage{amssymb}
\usepackage{array,color}
\usepackage{amsmath}
\usepackage{fullpage}

\usepackage{color}

\usepackage{mathrsfs}
\usepackage{amsmath,amssymb}
\usepackage{amsmath}
\usepackage{graphicx}
\usepackage{amsmath,amssymb,amsthm}
\usepackage{color}      
\usepackage{multirow}
\newtheorem{thm}{Theorem}[section]

\newtheorem{conj}[thm]{Conjecture}

\theoremstyle{definition}

\def \dpg { refined  defective painting game }

\begin{document}

\begin{frontmatter}

\title{Every planar graph is $1$-defective $(9,2)$-paintable}

\author[1]{Ming Han}
\address[1]{College of Mathematics, Physics and Information
Engineering\\ Zhejiang Normal University, Jinhua 321004, China.}

\author[1]{Xuding Zhu \fnref{label3}}
\fntext[label3]{ This research is supported by CNSF Grant 11571319. E-mail: xudingzhu@gmail.com.  }

\begin{abstract}
Assume $L$ is a $k$-list assignment of a graph $G$.
A $d$-defective $m$-fold $L$-colouring $\phi$ of   $G$ assigns to each vertex $v$ a set $\phi(v)$ of $m$ colours, so that $\phi(v) \subseteq L(v)$ for each vertex $v$,   and
 for each colour $i$, the set $\{v: i \in \phi(v)\}$ induces a subgraph of maximum degree at most $d$. 
 In this paper, we consider on-line list $d$-defective $m$-fold colouring of graphs, where the list assignment $L$ is given on-line, and the colouring is constructed on-line.
 To be precise, the $d$-defective $(k,m)$-painting game on  a graph $G$ is played by two players: Lister and Painter. Initially, each vertex has $k$ tokens and is uncoloured. In each round,
Lister chooses a set $M$ of  vertices  and removes one token from each chosen vertex.
Painter colours a subset  $X$  of $M$ which induces a subgraph $G[X]$ of maximum degree at most $d$.
A vertex $v$ is fully coloured if $v$ has received $m$ colours.
 Lister wins   if at the end of some round, there is a  vertex with no more tokens left and   is not fully coloured.
 Otherwise, at some round, all vertices are fully coloured  and
 Painter wins.
 We say $G$ is $d$-defective $(k,m)$-paintable if Painter has a winning strategy in this game.
 This paper proves that every planar graph is
 $1$-defective  $(9,2)$-paintable.
\end{abstract}


\end{frontmatter}


\section{Introduction}

A {\em $d$-defective $m$-fold colouring} of a graph $G$ is a colouring $\phi$ which assigns to each vertex
$v$ a set $\phi(v)$ of $m$ colours  so that
each colour class (i.e., each set of the form $\{v: i \in \phi(v)\}$ for some colour $i$)
induces a subgraph of maximum degree at most $d$. A $0$-defective $1$-fold colouring
of $G$ is simply a proper colouring of $G$,  a $0$-defective $m$-fold $k$-colouring of $G$ is called a $(k,m)$-colouring of $G$, and a $d$-defective $1$-fold colouring is also called
a {\em $d$-defective colouring } of $G$. A graph $G$ is called $d$-defective $k$-colourable if there
is a $d$-defective colouring of $G$ using $k$ colours.
Defective colouring of graphs was introduced by Cowen, Cowen and Woodall \cite{CCW1986}. They proved that
every outerplanar graph is $2$-defective $2$-colourable and every planar graph is $2$-defective $3$-colourable.

A {\em $k$-list assignment} is a mapping $L$ which assigns to each vertex $v$ a set
$L(v)$ of $k$ permissible colours. A {\em $d$-defective $m$-fold $L$-colouring} of $G$ is a $d$-defective
$m$-fold colouring $\phi$ of $G$ for which $\phi(v) \subseteq L(v)$ for every $v \in V(G)$.
A graph $G$ is {\em $d$-defective $(k,m)$-choosable} if for any $k$-list assignment $L$ of $G$,
there exists a $d$-defective $m$-fold $L$-colouring of $G$. If $G$ is $d$-defective $(k,1)$-choosable,
then we simply say that $G$ is {\em $d$-defective $k$-choosable}.
\u{S}krekovski  \cite{SKI1999} and Eaton
and Hull \cite{Eaton1999} independently extended the above result  to the list version and
proved that that every planar graph is  $2$-defective $3$-choosable and every outerplanar graph is $2$-defective $2$-choosable.
They both asked the question whether every planar graph is $1$-defective $4$-choosable.
One decade later, Cushing and Kierstead \cite{Kierstead2010} answered this question in the affirmative.

This paper studies  the on-line version of  defective multiple list colouring of graphs, defined through a two person game.

The {\em $d$-defective $(k,m)$-painting game on $G$} is played by two players: Lister and Painter.
Initially, each vertex $v$ has $k$ tokens and is uncoloured. In each round,
Lister chooses a set $M$ of  vertices   and removes one token from each chosen vertex.
Painter colours a subset  $X$  of $M$ which induces a subgraph $G[X]$ of maximum degree at most $d$.
A vertex $v$ is fully coloured if $v$ has received $m$ colours.
 Lister wins   if at the end of some round, there is a  vertex with no more tokens left and   is not fully coloured.
 Otherwise, at some round, all vertices are fully coloured  and
 Painter wins.
 We say $G$ is $d$-defective $(k,m)$-paintable if Painter has a winning strategy in this game.

 More generally,
 let $f: V(G) \to N$ be a mapping which assigns to each vertex $v$   a positive integer $f(v)$.
 The {\em $d$-defective $(f,m)$-painting game on $G$} is the same as the $d$-defective $(k,m)$-painting game, except that at the beginning of the game,
 each vertex $v$  has $f(v)$ tokens instead of $k$ tokens. If Painter has a winning strategy for the $d$-defective $(f,m)$-painting game on $G$, then we
 say $G$ is {\em $d$-defective $(f,m)$-paintbale}.
 If $m=1$, we simply say $G$ is {\em $d$-defective $f$-paintbale}.
The mapping $f$ is called a {\em token function}.


Observe that if $G$ is $d$-defective $(k,m)$-paintable, then $G$ is $d$-defective $(k,m)$-choosable.
Indeed, if $L$ is a $k$-list assignment of $G$, and in the $d$-defective $(k,m)$-painting game on $G$,
 Lister chooses  $M=\{v: i \in L(v)\}$ in the $i$th round, then  the colouring Painter constructed
 (using his winning strategy for the painting game) is a $d$-defective $(L,m)$-colouring of $G$.
On the other hand, there are graphs that are $d$-defective $(k,m)$-choosable but
 not $d$-defective $(k,m)$-paintable. For example,
 it is proved in  \cite{Eaton1999,SKI1999} that every planar graph is $2$-defective $3$-choosable,
 however, it is recently shown in \cite{GGHZ2016} that there are planar graphs that are not
$2$-defective $3$-paintable.

A graph is  $(f,m)$-choosable if it is   $0$-defective   $(f,m)$-choosable, 
and  $(f,m)$-paintable if it is   $0$-defective   $(f,m)$-paintable,
A  graph is $f$-choosable if it is $(f,1)$-choosable and $f$-paintable if it is $(f,1)$-paintable.

Painting game on planar graphs and local planar graphs have been studied a lot in the literature.
It is known that every planar graph is $5$-paintable \cite{Schauz2009}
and $3$-defective $3$-paintable \cite{GGHZ2016}, every   planar graph of girth at least $5$ is $3$-paintable
 \cite{CZ},  every outerplanar graph is $2$-defective $2$-paintable \cite{HZ2015}, and for any surface $S$, there is a constant $w$ such that every graph embedded in $S$ with edge-width at least $w$ is
 $5$-paintable \cite{HZ2014} and  $2$-defective $4$-paintable \cite{HZ2015}.

Multiple list colouring and painting game on graphs have been studied in a few papers. It was proved by Alon, Tuza and Voigt  \cite{ATV97} that if $G$ is $(a,b)$-colourable, then for some integer $m$, $G$ is $(am,bm)$-choosable. In particular, for each planar graph $G$, there is an integer $m$ such that $G$ is $(4m,m)$-choosable. However, the 
integer $m$ depends on $G$. It is recently proved in \cite{Zhu2016} that for any positive integer $m$, there is a planar graph $G$ which is not $(4m,m)$-choosable.
The following conjecture is posed in \cite{Zhu2016}:
 
\begin{conj}
	\label{conj1}
   There is a constant integer $m$ such that every planar graph $G$ is $(5m-1,m)$-choosable.
	\end{conj}	

If Conjecture \ref{conj1} is true, then the constant $m$ needs to be at least $2$. The following conjecture (also posed in \cite{Zhu2016}) says that $2$ is enough.   
	
	\begin{conj}
		\label{conj2}
		Every planar graph $G$ is $(9,2)$-choosable.
	\end{conj}
	 

For painting game, it was proved by Gutowski \cite{Gutowski2011} that if $G$ is $(a,b)$-colourable, then for any $\epsilon > 0$, there is an integer $m$ such that $G$ is $((a+\epsilon)m,bm)$-paintable. Thomassen's proof \cite{Tho1994} that every planar graph is $5$-choosable can be easily adopted to prove that for any positive integer $m$, every planar graph is $(5m, m)$-paintable. Analog to Conjecture \ref{conj2}, we have the following conjecture:

\begin{conj}
	\label{conj3}
Every planar graph $G$ is $(9,2)$-paintable.
\end{conj}

Conjecture  \ref{conj3} implies  Conjecture \ref{conj2}, which in turn is stronger than Conjecture \ref{conj1}. It follows from the Four Colour Theorem  that for any integer $m$, every planar graph is $(4m,m)$-colourable. It was proved in \cite{HRS1973} (without using the Four Colour Theorem) that every planar graph $G$ is $(5m-1,m)$-colourable with $m = |V(G)|+1$.  However, there is no direct proof  of the fact that every planar graph is $(5m-1,m)$-colourable for some constant $m$.
This suggests that these  conjectures might be difficult. 

In this paper, we prove the following result, which can probably be viewed as a weak support of Conjecture \ref{conj3}.  

\begin{thm}\label{key}
Every planar graph  is $1$-defective $(9,2)$-paintable.
\end{thm}

One natural question suggested by the result of Cushing and Kierstead is that whether every planar graph is $1$-defective $(8,2)$-choosable. A long standing difficult conjecture of Erd\H{o}s-Rubin-Taylor states that if a graph $G$ is $k$-choosable, then for any positive integer $m$, $G$ is $(km,m)$-choosable \cite{ERT79}. Analog to this conjecture,  one may conjecture that if a  graph $G$ is $d$-defective $k$-choosable, then for any integer $m$, $G$ is $d$-defective $(km,m)$-choosable.    Erd\H{o}s-Rubin-Taylor's conjecture has received a lot of attention, however, little progress has been made on this conjecture in the past 40 years. Approach to the analog conjecture about defective list colouring seems to meet the same difficulties. Our result implies that  
every planar graph is  $1$-defective $(9,2)$-choosable, which is  weaker than the conjectured result.  
Another question suggested by the result of Cushing and Kierstead is that whether every planar graph is $1$-defective $4$-paintable. This question is still open. The best result in this direction is that for any surface $S$ there is a constant $w$ such that any graph embedded in $S$ with edgewidth at least $w$ is   $2$-defective $4$-paintable \cite{HZ2015}.

\section{Preliminaries}

Let $C$ be the boundary of a plane graph. For $x, y \in C$,  denote by $C[x,y]$  the path on $C$ from $x$ to $y$ in the clockwise direction. Let $C(x,y)=C[x,y]-\{x,y\}$, $C(x,y]=C[x,y]-\{x\}$ and $C[x,y)=C[x,y]-\{y\}$.
For a cycle $C'$ of $G$,  denote by ${\rm int}[C']$  the subgraph of $G$ induced by all vertices   inside or on $C'$, and denote by ${\rm int}(C')$ the subgraph of $G$
induced by vertices inside $C'$ (but not on $C'$). For a vertex $v$ of $G$, let $N_G(v)$ be the set of neighbours of $v$,
and for a subset $A$ of $V(G)$, let $N_G(A)=\cup_{v \in A}N_G(v)$. If the graph $G$ is clear from the text, we write $N(v)$ and $N(A)$ for $N_G(v)$ and $N_G(A)$, respectively.
If $P$ and $P'$ are two paths of $G$ that are vertex disjoint except that the last vertex of $P$ is the same as the first vertex of $P'$, then the concatenation of $P$ and $P'$, written as $P \cup P'$,
is the path which is the union of $P$ and $P'$.
If $P$ and $P'$ are vertex-disjoint, except that the last vertex of $P$ is the   first vertex of $P'$ and
the last vertex of $P'$ is  the first vertex of $P$, then $P \cup P'$ is a cycle.

For the remainder of this paper, we restrict to $1$-defective $(9,2)$-painting game.
For the purpose of using induction,
instead of proving   Theorem \ref{key} directly, we prove
a stronger and more complicated result.

 If a vertex $v$ is coloured in a certain round, then its defect in that round is the number of  neighbours coloured in this round.
If two adjacent vertices $u$ of $v$ are  coloured   in the same round, then $u$ is said to  {\em  contribute one defect} to $v$ and also 
{\em receive  one defect} from $v$.

In the $1$-defective $(9,2)$-painting game,   each vertex $v$ has $9$ tokens,
and in a round $v$ is coloured, it can have defect at most $1$.
 We shall consider a   refined  defective painting game.

One feature of the refined defective painting game is that tokens assigned to vertices carry values.
A token can be a {\em one-dollar token} or a {\em two-dollar token}.
When Lister takes away from $v$ an $s$-dollar token,
Painter is allowed to colour $v$ in this round, under the restriction that it can have defects at most   $s-1$.

We need this feature, because in playing the game, $V(G)$ is divided into several parts, possibly with intersecting boundaries, and
for each part, the game is played on a smaller graph and hence Painter has a winnings strategy by induction hypothesis.
A vertex $v$   coloured in a round may receive one defect from  one  part, and hence  cannot receive any defect from the other parts. 
So for our stronger technical result to be proved by induction,
vertices on the boundary may have one-dollar tokens. Moreover, vertices on the boundary may have fewer tokens due to its neighbours coloured
in the other parts. In the refined defective painting game, the token function  $f$ assigns to each vertex $v$ of $G$ a pair of integers $(a,b)$, which means that   $v$ is assigned $a$ one-dollar tokens and $b$ two-dollar tokens. If $f(v)=(a,b)$, then $v$ is called an $(a,b)$-vertex.

Assume $G$ is a plane graph and $C=c_1c_2\cdots c_n$ is the boundary of the infinite face of $G$, and $f$ is a token function on $G$.
We say the
  token function $f$   is {\em valid for the refined defective painting game on $G$ with special vertices $c_1,c_2,c_n$} (shortened as
   {\em valid for   $(G;c_1,c_2,c_n)$}) if the following hold:
\begin{enumerate}
\item[(F1)] $f(c_1)=f(c_2)=(2,0)$.
\item[(F2)] If $N_G(c_n) \cap N_G(c_1) \cap (V(C)-\{c_2\}) = \emptyset$, then $f(c_n)=(0,2)$. Otherwise, $f(c_n)=(4,0)$.
\item[(F3)] $f(c_i) = (2,5)$ for $3 \le i \le n-1$,
\item[(F4)] $f(v) = (0,9)$ for every interior vertex $v$.
\end{enumerate}

Let $f$ be a valid token function of $G$.
The {\em  refined defective painting game on $G$ with special vertices $c_1,c_2,c_n$}
is played by two players: Lister and Painter.

Initially, for  each  vertex $v$ of $G$, if $f(v)=(a,b)$, then $v$ is   assigned
 $a$ one-dollar tokens and $b$ two-dollar tokens. For each vertex $v$, we  denote by
 $\beta(v)$   the number of rounds in which $v$ is coloured, and let $\beta(v)=0$ at the beginning of the game.

In each round, Lister's move is to present a mapping $M: V \to \{0,1,2\}$, and if $M(v) =j >0$, then remove one $j$-dollar token
from $v$,
Lister's move is legal if
\begin{enumerate}
 \item[(L1)] $M(c_n)M(c_1)=0$.
 \item[(L2)] If $c_n$ and $c_2$ are adjacent, then $M(c_n)M(c_2)=0$.
 \item[(L3)] If $M(v)=j$, then before this round, $v$ has at least one $j$-dollar token.
 \end{enumerate}
If $M(v) =s \ge 1$, then we say $v$ is {\em marked with an $s$-dollar token} in this round.
The set of vertices $v$ with $M(v) \ge 1$ are {\em marked vertices} in this round. Condition  (1)   implies that
  $c_n$ and $c_1$ are never marked in the same round.
 Condition  (2)   implies that if $c_2$ and $c_n$ are adjacent, then $c_n$ and $c_2$ are also never marked in the same round.
We denote by  $e^*$   the edge $c_1c_2$.
 Painter's move is to chooses a subset $X$ of marked vertices, and for each $v \in X$, increase $\beta(v)$ by $1$.
 Painter's move is legal if  $$\forall v \in X, |X \cap N_{G-e^*}(v)| \le M(v)-1. \eqno(P1)$$

Note that  (P1) implies that if $M(v)=0$, then
$v \notin X$. The inequality (P1) bounds the defect of each coloured vertex $v$.
However, we only count the defect in the graph $G-e^*$.
It may happen that $c_1$ and $c_2$ are coloured in the same round. In this case,
the defect contributed by $c_1$ to $c_2$ and by $c_2$ to $c_1$ are not counted in (P1). Since $c_1$ and $c_2$ are $(2,0)$-vertices,
no other neighbour contributes defect to them. So even if $c_1$ and $c_2$ are coloured in the same round, they have
defect only $1$ in $G$.

A vertex $v$ is called {\em fully coloured} if $\beta(v) \ge 2$.
If at the end of some round there is a vertex with no tokens left and is not fully coloured,
then Lister wins the game. Otherwise, at the end of some round,
all vertices are fully coloured  and Painter wins the game.

For two pairs $(a,b)$ and $(a',b')$ of non-negative integers, we write $(a',b') \preceq (a,b)$ if $b' \le b$ and $a'+b' \le a+b$. For example, $(2,5) \preceq (0,7) \preceq (0,9)$.
For $f, f': V(G) \to N^2$, we write $f' \preceq f$ if for any vertex $v$, $f'(v) \preceq f(v)$.
In our proofs, we shall frequently need to show that a token function $f$ is valid for $(G; c_1, c_2, c_n)$.
For convenience, to prove $f$ is valid for $(G; c_1, c_2, c_n)$,  it suffices to show that $f(v)\succeq (2,0)$ for $v\in \{c_1,c_2\}$,
 $f(c_n)\succeq(0,2)$ or $f(c_n)\succeq(4,0)$,
$f(v)\succeq (2, 5)$ for $v \in V(C)-\{c_1,c_2,c_n\}$ and
$f(v)\succeq(0,9)$ for $v\in V(G)-V(C)$, as we may drop off extra tokens or devalue some tokens if needed.
(To devalue the token of vertex $v$ in a round means to change $M(v)$ to $\max\{0, M(v)-1\}$).

\begin{thm}\label{main}
Let $G$ be a plane connected graph with boundary walk $C$. If $f$ is a valid token function for $(G;c_1,c_2,c_n)$, then Painter has a winning strategy for the\dpg on $G$ with special vertices $c_1,c_2,c_n$.
\end{thm}

It is obvious that Theorem \ref{main} implies Theorem \ref{key}.
 
 The proof of Theorem \ref{main} is by induction on the number of vertices of $G$.
 Depending on the structure of $G$, we shall consider a few cases.
 In each case, we shall consider some induced subgraphs $G_1, \ldots, G_k$ of $G$ such that $V(G_1) \cup \ldots \cup V(G_k) =V(G)$
 (these graphs need not be vertex disjoint, they may intersect on their boundaries).
 We shall first play a painting game on $G_1$, then on $G_2, G_3$ and so on.
 In a given round, assume Lister has made  a move $M$. We denote by $M_i$ Lister's move on $G_i$,
 and denote by $X_i$ Painter's move on $G_i$. The union $X = \cup_{i=1}^kX_i$ is painter's move on $G$ in
 this round. 
 
 The move $M_i$ depends on $M$ and may also depends on $X_1, X_2, \ldots, X_{i-1}$. We shall always have $M_i(v) \le M(v)$ for all $v \in V(G_i)$.
 It may happen that $M_i(v) < M(v)$,  because $v$ may have  a neighbour in another subgraph $G_j$ which is coloured
 in this round. 

Each $G_i$ has fewer number of vertices than $G$. By induction hypothesis,
 Theorem \ref{main} holds for $G_i$.   To play the game on $G_i$   by induction, we shall choose three consecutive vertices,
 say $c_{i,1}, c_{i,2}, c_{i,n_i}$,
 on the boundary $C_i$ of $G_i$ as
 special vertices.
 Lister's moves $M_i$ for the games on $G_i$ should be legal and should guarantee that
 enough tokens are assigned to vertices of $G_i$ so that induction can be applied.
 Assume $v$ is a vertex of $G_i$. Let $f_i(v)=(a,b)$ if
 in the game on $G_i$,  $v$ is marked with one-dollar token in $a$ rounds, and marked with two-dollar token in $b$ rounds.
 With an abuse of notation, we call $f_i$ the token function of $G_i$.
 (In the definition of a token function $f$, $f(v)$ should be the tokens assigned to $v$ before the play of the game.
 Here $f_i(v)$ is the tokens actually used by Lister on the game on $G_i$.)

 To use induction on $G_i$, we need to show that $f_i$ is valid for $(G_i;c_{i,1},c_{i,2},c_{i,n_i})$.
 We say $f_i$ is {\em valid at $v$} (for the game on $G_i$) if the set of tokens assigned to $v$ is 
 enough, i.e.,  $f_i(v)  \succeq (2,0)$ for $v\in \{c_{i,1},c_{i,2}\}$,
 $f_i(c_{i,n_i})\succeq(0,2)$ or $f_i(c_{i,n_i})\succeq(4,0)$,
$f_i(v)\succeq (2, 5)$ for other boundary vertex $v$ of $G_i$ and
$f_i(v)\succeq(0,9)$ for interior vertices $v$ of $G_i$. To prove that $f_i$ is valid, we need to   show that $f_i$ is valid at every vertex
 $v \in V(G_i)$.

 There are cases in which some $G_i$ has very simple structure (for example, $G_i$ may be a single vertex). Instead of playing the painting game on $G_i$ by induction hypothesis, we
   colour the vertices of $G_i$ by some explicitly described simple rules.

On deciding Lister's move and Painter's move on each $G_i$, we
 also need to make sure that Painter's moves on $G_i$'s are consistent (i.e., if a vertex $v$ belongs to two subgraphs $G_i$ and $G_j$, then
it should be either coloured in both games or not coloured in both games), and the move on $G$ must be legal.

\section{Proof of Theorem \ref{main}}

Based on the structure of $G$, we divide the proof into a few cases.  If $G$ has a vertex $v$ of degree $1$, 
then   Painter's winning strategy for the\dpg  on $G-v$ can be easily adapted to a winning strategy for the game   on $G$.   
So we assume $G$ had minimum degree $\delta(G) \ge 2$. 

We assume that in the process of the game, if a vertex $v$ is fully coloured (i.e., coloured in two rounds), then in later moves, $v$ will not be marked.  In the games on subgraphs $G_i$, we shall define Lister's move as $M_i$ by certain formulas. We make the convention that if $v$ is fully coloured before this move, then $M_i(v)$ is set to $0$, no matter what the formula says. This is to forbid a vertex to be coloured in more than two rounds.


\bigskip
\noindent
{\bf Case 1.} $G$ has a cut-vertex.

\begin{proof}
Assume $G$ contains a cut-vertex $w$. 
Let $G_1, G_2$ be subgraphs of $G$ with $V(G_1) \cup V(G_2)=V(G)$ and
$V(G_1) \cap V(G_2)=\{w\}$. 

First we consider the case that all of  $c_1,c_2,c_n $ are contained in one of the subgraphs $G_1$ and $G_2$.
We assume $c_1,c_2,c_n  \in V(G_1)$, as depicted in Figure 1(a).  As $\delta(G) \ge 2$, we can choose the cut-vertex $w$ so that $w$
has at least two neighbours in $G_2$.  Let $w',w''$ be the neighbours of $w$ on the boundary of $G_2$. 

\begin{figure}[!ht]
  \begin{center}
    \includegraphics[scale=0.45]{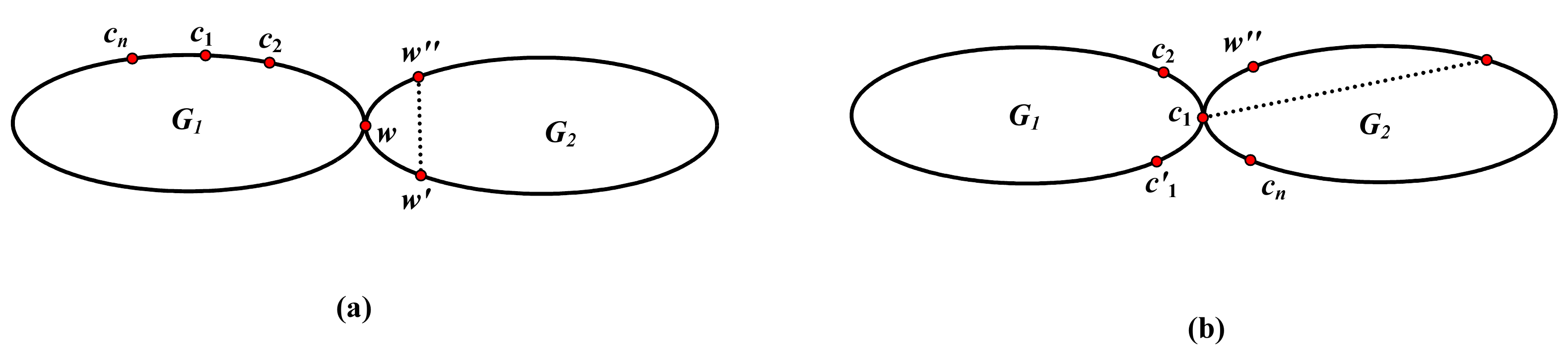}
\caption{$G$ has a cut-vertex  $w$.}

   \end{center}
\end{figure}

Assume Lister has made a move $M$.
In the first phase,  Painter plays the\dpg  on $G_1$ with special vertices $c_1,c_2,c_n$, where Lister's move on $G_1$ is $M_1=M|_{V(G_1)}$. 
So $f_1=f|_{V(G_1)}$ and hence is valid for  $(G_1;c_1,c_2,c_n)$.
Thus, Painter  has the winning strategy for the game on $G_1$ by induction hypothesis and can choose a subset $X_1$ based on the winning strategy.

In the second phase, Painter plays the\dpg on   $G_2$
with special vertices $w,w'',w'$,  where 
  Lister's move $M_2$ is the restriction of $M$ to $V(G_2)$, except that
  \begin{itemize}
  	\item  $M_2(w)=\begin{cases}
  	1&\text{if}~w \in X_1,\\
  	0&\text{otherwise}.
  	\end{cases}$
  	\item For $w',w''$,
  	\subitem if $w'w''\notin E(G)$, then for $v \in \{w', w''\}$, $M_2(v)=0$ whenever 
  	$w\in X_1$.
  	\subitem if $w'w''\in E(G)$, then 
    $M_2(w')=\begin{cases}
  	      2&\text{if}~M(w')=2~\text{and}~w\notin X_1,\\
  	      0&\text{otherweise},
  	      \end{cases}$\\
  	      and $ M_2(w'')=0$ whenever $M_2(w')\ne 0$ or $w\in X_1$.
  \end{itemize}


In the game of $G_2$, $w,w'',w'$ play the roles of $c_1,c_2,c_n$, respectively.

 Recall that for $v \in V(G_2)$, $f_2(v)=(a,b)$ means that for the game on $G_2$, Lister marks $v$ with one-dollar token
 in $a$ rounds and marks $v$ with two-dollar token in $b$ rounds. We need to show that $f_2$ is valid for $(G_2;w,w'',w')$.
 It is obvious that $f_2$ is valid at $v$ for any vertex $v \in V(G_2)-\{w,w',w''\}$. 
 
By definition of $M_2$,   $w'$ and $w$ (and in case  $w'w''\in E(G)$, then $w'$ and $w''$) are not marked together in the same round, so for the game on $G_2$, Rules (L1) and (L2) are satisfied.

Also by definition of $M_2$,  $w$ is marked with a one-dollar token if and only if $w$ is coloured in the game on $G_1$.
As $w$ is coloured in two rounds, we know that $f_2(w) = (2,0)$, hence $f_2$ is valid at $w$.
 
Assume $w'w'' \in E(G)$.  It follows from the definition of $M_2$ that
 in a certain round, if $w$ is coloured   and $w'$ is marked in the game on $G$, then $w'$ is not marked
 in the game on $G_2$. In this case, we say $w'$ lost one token due to $w$. As $w$ is coloured in two rounds, $w'$ may lose two tokens due to 
 $w$. Also by the definition of $M_2$, we know that  if $w'$ is marked with a one-dollar token in the game on $G$, then it is not marked 
 in the game on $G_2$. So $w'$ lose two one-dollar tokens.  As $w'$ is a $(2,5)$-vertex in the game on $G$, $f_2(w') \succeq (0,3)$. 
 Since $w'$ plays the role of $c_n$ in $G_2$, and  $w'w'' \in E(G)$   implies that
  $w'$ and $w$ have no common neighbours other than $w'$ on the boundary of $G_2$,    $f_2$ is valid at $w'$.
 Similarly, it follows from the definition of $M_2$ 
 that $w''$ may lose four tokens due to $w$ and $w'$, and hence $f_2(w'') \succeq (2,1)$ and $f_2$ is valid at $w''$.  

Assume $w'w''\notin E(G)$. Then  each vertex of $w',w''$ may lose two tokens due to $w$, and   $f_2(w'), f_2(w'') \succeq(2,3)$, and hence
 $f_2$ is also valid at $w'$ and $w''$.

In any case, $f_2$ is valid for $(G_2; w,w'',w')$.
By induction hypothesis, Painter has the winning strategy for the game on $G_2$ and can choose a subset $X_2$ based on the winning strategy.
 The union $X=X_1\cup X_2$ is Painter's move for  the game on $G$.

Note that $w$ is coloured in the game on $G_2$ if  and only if $w$ is coloured in the game on $G_1$. So Painter's moves $X_1$  and $X_2$ are consistent.
Moreover, as $w,w''$ play  the roles of $c_1,c_2$ in the game on $G_2$ and $w''$ is not marked in the same round as $w$,  $w$ receives no defect from $G_2$. Hence $w$ satisfies Rule (P1). Therefore,  Painter's move for the game on $G$  is legal.

Next we consider the case that none of $G_1$ and $G_2$ contains all the vertices $c_1, c_2, c_n$.
In this case, we have $w=c_1$. Assume $c_2$ is contained in $G_1$ and $c_n=w'$ is contained in $G_2$, as depicted  in Figure 1(b).
We may assume that $w''$ (as shown in Figure 1(b)) and $c_n$ are distinct, as otherwise we could choose $c_n$ as the cut-vertex. 

For $i=1,2$, let $c'_i$ be the other neighbour of $c_i$ on the boundary of $G_1$. Similarly, $c'_1$ and $c_2$ are distinct
vertices, as otherwise we could choose $c_2$ as the cut-vertex. We know that either $c'_1$ is not adjacent to $c_2$ or $c'_2$ is not 
adjacent to $c_1$. Without loss of generality, we may assume that $c'_1$ is not adjacent to $c_2$.

Assume Lister has made a move $M$.
In the first phase,  Painter plays the\dpg on   $G_1$
with special vertices $c_1,c_2,c'_1$, where Lister's move $M_1$ is the restriction of $M$ on $V(G_1)$, except that  $M_1(c'_1)=0$ whenever $M(c_1)\ne0$. 

 By the definition of $M_1$, $M_1(c'_1)M_1(c_1)=0$ and $c_2c'_1\notin E(G)$ (by the choice of $c'_1$), so Lister's move is legal. 
  The vertex $c'_1$  may lose two tokens due to $c_1$, and hence $f_1(c'_1)\succeq(2,3)$ and $f_1$ is valid at $c'_1$. Obviously, $f_1$ is valid at other vertices of $G_1$. Let $X_1$ be the subset chosen by Painter based on his winning strategy for the game on $G_1$.

In the second phase,  Painter plays the\dpg on   $G_2$
with special vertices $c_1,w'',c_n$.
In the game on $G_2$,   Lister's move   $M_2$ is the restriction of $M$ to $V(G_2)$, except  that  $M_2(w'')=0$ whenever $M(c_1)\ne 0$ or $w''c_n\in E(G)$ and $M(c_n)\ne 0$. 
As $c_n$ and $c_1$ are not  marked in the same round (and not marked in the same round as 
$w''$  in case $w''c_n\in E(G)$),  Lister's move $M_2$ is legal. If $w''c_n\in E(G)$, then $f(c_n)=(0,2)$ since $N(c_n)\cap N(c_1)\cap (V(C)-\{c_2\})=\emptyset$. Thus, the vertex   $w''$   loses at most four tokens,  and hence $f_2(w'')\succeq(2,1)$. If $w''c_n\notin E(G)$, then $w''$   may lose two tokens due to $w$, and hence $f_2(w'')\succeq(2,3)$. It is obvious that $f_2$ is valid at other vertices of $G_2$. 
Hence Painter has the winning strategy for the game on $G_2$ by induction hypothesis and could choose a subset $X_2$ in this round.

 Painter's moves $X_1$  and $X_2$ are consistent since $c_1=w$ has just two tokens and would be coloured when it is marked.  Moreover, as $c_1,w''$ play the roles of $c_1,c_2$ in the game on $G_2$ and $w''$ is not marked in the same round as $c_1$ by the definition of $M_2$,  $c_1$ receives no defect from $G_2$. Meanwhile, $c_1,c_2$ play the roles of $c_1,c_2$ in the game on $G_1$.  So $c_1$ receives no defect from $G-\{e^*\}$. Therefore,  Painter's move  $X=X_1\cup X_2$ for the game on $G$  is legal.

\end{proof}

Thus $G$ is $2$-connected and the boundary $C$ of $G$ is a cycle.
We may assume that $G$ is near-triangulated, i.e., each facial cycle of $G$ other than $C$  is a triangle.

\bigskip
\noindent
{\bf Case 2.} $C$ is a triangle.

\begin{proof}
If $C$ is a triangle $(c_1,c_2,c_3)$, then let $G_1=\{c_3\}$ and $G_2=G-\{c_3\}$ (as depicted  in Figure 2).
\begin{figure}[!ht]
  \begin{center}
    \includegraphics[scale=0.45]{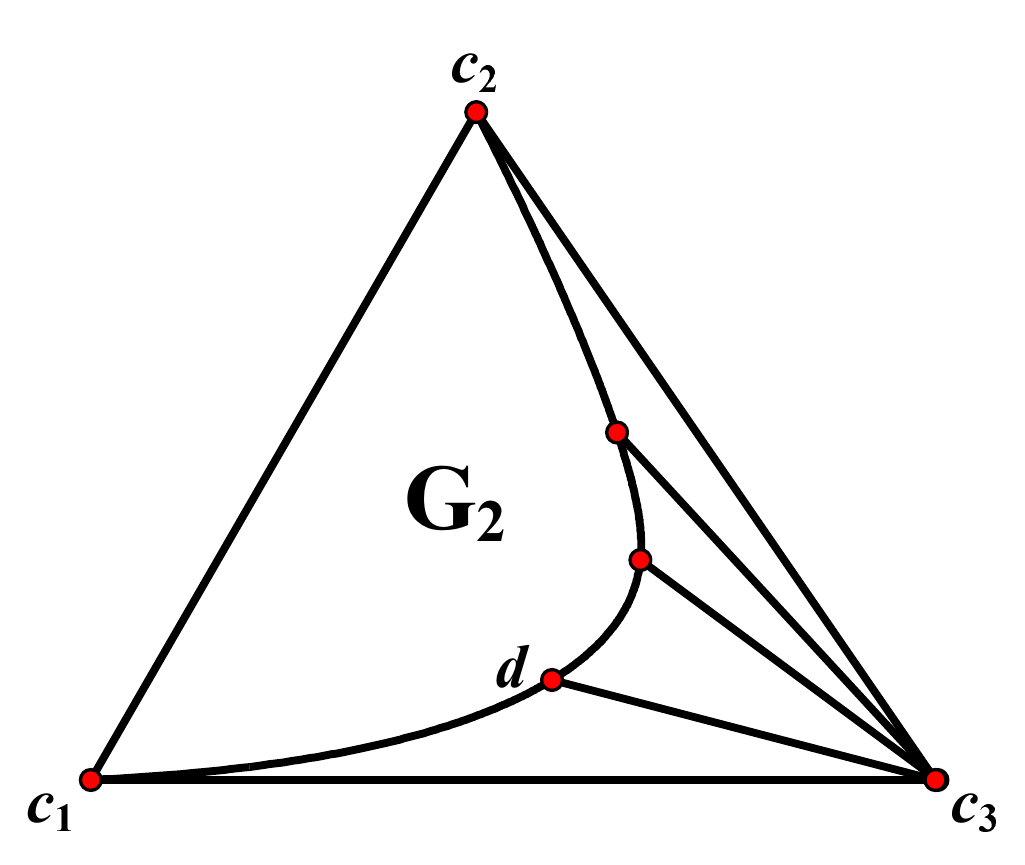}
\caption{The boundary of $G$ is a triangle.}

   \end{center}
\end{figure}
 Note that $N(c_1)\cap N(c_3)\cap (V(C)-\{c_2\})=\emptyset$, so $f(c_3)=(0,2)$.
Assume Lister has   made a move $M$.  
In the first phase, Painter plays the game on $G_1$ which contains just one vertex $c_3$ and colours $c_3$ if and only if $M(c_3)=2$.
In the second phase, Painter plays the\dpg on $G_2$ with    special vertices $c_1,c_2,d$, where $d$ is the other neighbour of $c_1$ on the boundary of $G_2$.
Lister's move $M_2$ is the restriction of $M$ to $V(G_2)$, except that  for $v \in N_G(c_3)-\{c_1,c_2,d\}$, $M_2(v)=0$ whenever $c_3 \in X_1$ and 
$M_2(d)=0$  whenever $M(c_1)+M(c_3)\ne 0$ or $c_2d\in E(G)$ and $\sum_{i=1}^3M(c_i)\ne 0$.

Thus, the vertex $d$ is not marked in the same round as $c_1$ (and not marked in the same round as 
$c_2$  in case $dc_2\in E(G)$). So Lister's move for the game on $G_2$ is legal.  Moreover, if $dc_2\in E(G)$, then $d$ may lose six  tokens  due to $c_1,c_2$ and $c_3$, and hence $f_2(d)\succeq(0,3)$.   As $dc_2\in E(G)$, $N(c_1)\cap N(d) \cap (V(C_2)-\{c_2\})=\emptyset$. Thus, the token function $f_2$ is valid at $d$.
If $dc_2\notin E(G)$, then  $d$ may lose four tokens due to $c_1$ and $c_3$, and $f_2(d)\succeq(0,5)$.
So $f_2$ is also valid at $d$.
By induction hypothesis, Painter can choose a subset $X_2$ in this round.
  
Since $V(G_1)\cap V(G_2)=\emptyset$, so Painter's moves on $G_1$ and $G_2$ are consistent. As Lister's move in the game on $G$ is legal,  $c_3$ is never marked in the
same round as $c_1$ and $c_2$. By the definition of $M_2$, $c_3$ is never marked in the same round as any other neighbour.
Therefore,  $c_3$ contributes  no defect to   $G_1$ and also receives no defect from $G_1$,
and  Painter's move  $X=X_1 \cup X_2$ is legal for the game on $G$.

\end{proof}

\bigskip
\noindent
{\bf Case 3.} $G$  has a separating triangle.

\begin{proof}
Assume $C'=(x,y,z)$ is a separating triangle of $G$.
Let $G_1=G-{\rm int}(C')$, $G_2={\rm int}[C']$.

Assume Lister has made a move $M$.
In the first phase, Painter plays the\dpg  on $G_1$ with  special vertices $c_1,c_2,c_n$, where 
Lister's move $M_1=M|_{V(G_1)}$.  Lister's move is obviously legal, and the token function $f_1=f|_{V(G_1)}$ is valid.
Let $X_1$ be the subset chosen by Painter based on  the  winning strategy for the game on $G_1$.

In the second phase,   
Painter  plays the game  on $G_2$ in the same way as the game   in Case 2, with $x,y,z$ play the roles of $c_1,c_2,c_3$, respectively.
We mark $x,y,z$  in the game on $G_2$ if and only if they are coloured in $G_1$. So Painter's moves on $G_1$ and $G_2$ are consistent.  
By the proof of Case 2, we know that $x, y,z$ receive no defect from $G_2-\{x,y,z\}$. So Painter's move is legal. 
\end{proof}

\bigskip
\noindent
{\bf Case 4.} $G$ has  a chord.

\begin{proof}
We consider the following two subcases.

\bigskip
\textbf{Case 4.1.}: $G$ has a chord which is not incident to $c_1$.

We choose a chord $c_ic_j$ ($2\le i<j\le n$) so that  $G_2$ is chordless, where  $G_1={\rm int}[C[c_j,c_i]\cup \{c_ic_j\}]$
and $G_2={\rm int}[C[c_i,c_j]\cup \{c_jc_i\}]$  (as depicted in Figure 3(a)).

\begin{figure}[ht]
  \begin{center}
    \includegraphics[scale=0.45]{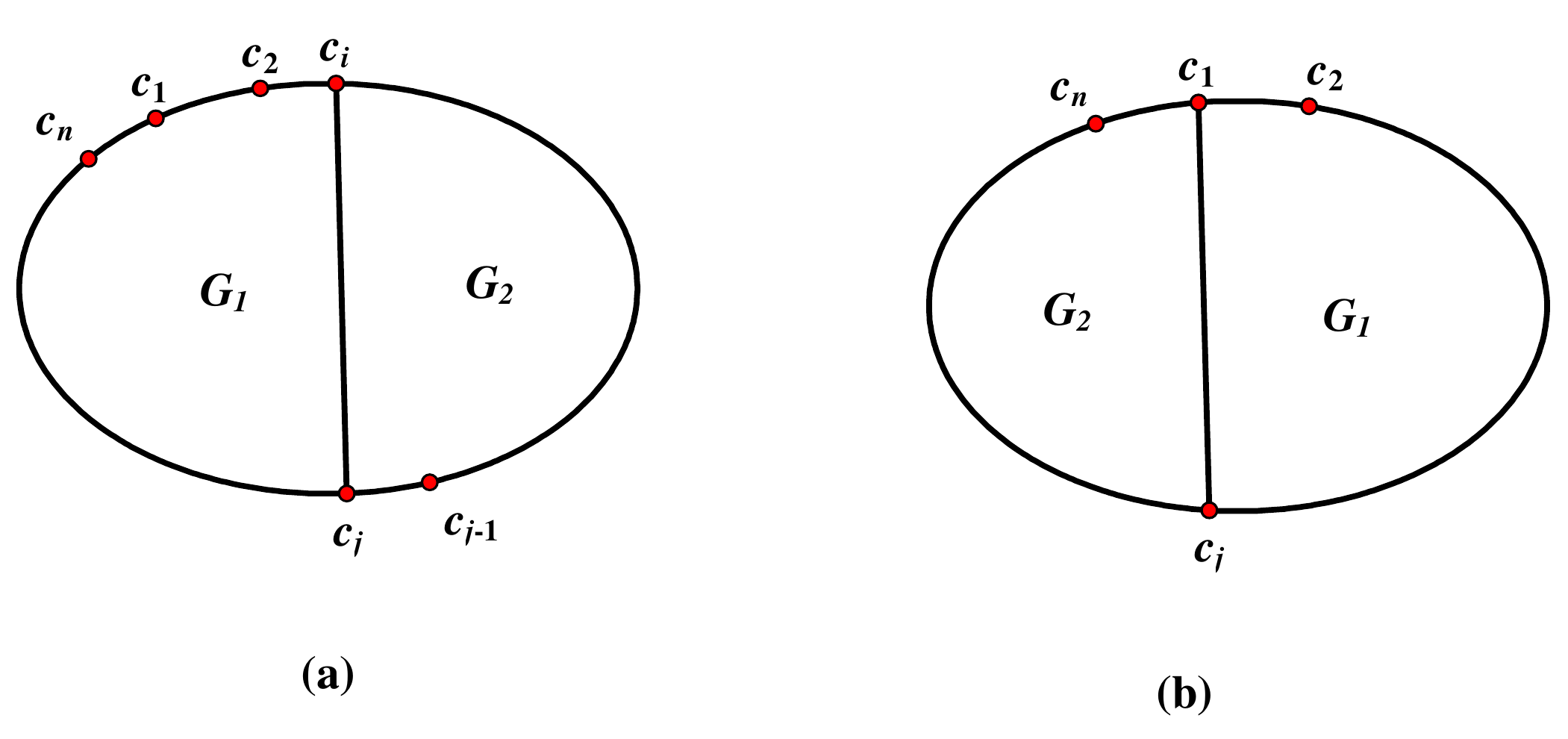}
\caption{$G$ has a chord $c_ic_j$.}

   \end{center}
\end{figure}

 Assume Lister has made  a move $M$.
In the first phase, Painter plays the\dpg on $G_1$ with special vertices $c_1,c_2,c_n$, where   Lister's move $M_1=M|_{V(G_1)}$. 
Let $X_1$ be the subset chosen by Painter  according to his winning strategy for the game on $G_1$.
In the second phase, if the boundary of $G_2$ is a triangle (i.e., $j=i+2$), then the interior of $G_2$ is empty as  $G$ has no separating triangle. Painter colours $c_{j-1}$ if and only if $M(c_{j-1})\ge 1$ and $c_i,c_j$ are not coloured. So $c_{j-1}$  loses at most four tokens  due to $c_i,c_j$ and hence can be fully coloured.
If $G_2$ is not a triangle, then Painter plays the\dpg on $G_2$ with special vertices $c_i,c_j,c_{j-1}$, where   Lister's move  $M_2$ is the restriction of $M$ to $V(G_2)$, except that  $M_2(c_{j-1})=0$ whenever $c_j \in X_1$ 
and for $v \in \{c_i,c_j\}$, 
$$M_2(v)=\begin{cases}
		1&\text{if}~v\in X_1,\\
		0&\text{otherwise}.
		\end{cases}
		$$
		
	By the definition of $M_2$, $c_{j-1}$ is not marked in the same round as $c_j$.  Since $G_2$ is chordless, the vertex $c_{j-1}$ is not adjacent to $c_i$. Thus Lister's move is legal. Moreover, the vertex $c_{j-1}$ may lose two tokens, and hence  $f_2(c_{j-1})\succeq(2,3)$. For $v\in\{c_i,c_j\}$, it would be fully coloured in the first phase, and hence $f_2(v)=(2,0)$. The other vertices of $G_2$ obviously have enough tokens. So the token function $f_2$ is valid and Painter can choose a subset $X_2$ by induction hypothesis.
	
	The common vertices $c_i$ and $c_j$ of $G_1$ and $G_2$ are coloured consistently in both games. 
	Moreover,  as $c_i$ and $c_j$  play  the roles of $c_1$ and $c_2$, they receive no defect from   $V(G_2)-\{c_i,c_j\}$. 
	Therefore, Painter's move for the game on $G$ is legal.

\bigskip
\textbf{Case 4.2}: Every chord of $G$ is incident to $c_1$.

 Assume $c_1c_j$ is a chord of $G$.  Let   $G_1={\rm int} [C[c_1,c_j]\cup \{c_jc_1\}]$
 and $G_2={\rm int}[C[c_j,c_1]\cup \{c_1c_j\}]$ (as depicted  in Figure 3(b)).

 Assume Lister has made a move $M$.

 In the first phase,  if $G_1$ is a triangle,  then Painter colours $c_1,c_2$ when  they are marked and  colours $c_j$ when $M(c_j)\ge 1$ and   $M(c_1)+M(c_2)=0$.  As $c_j$ has seven tokens and may lose at most four tokens, it can be fully coloured.
Assume $G_1$ is not a triangle.  Painter plays the\dpg on $G_1$ with special vertices $c_1,c_2,c_j$,  where  Lister's move $M_1$ is the restriction of $M$ to $V(G_1)$, except that $M_1(c_j)=0$ whenever $M(c_1) \ne 0$.

By the definition of $M_1$, $c_j$  is not marked  in the same round as  $c_1$.  Since $G_1$ is not a triangle (and hence its boundary is not a triangle)
and every chord of $G$ is incident to $c_1$, $c_j$ is not adjacent to $c_2$. So  Lister's move for the game on $G_1$ is legal. Moreover, the vertex $c_j$   may lose two tokens, and  hence $f_1(c_j)\succeq(2,3)$.
The other vertices of $G_1$ obviously have enough tokens. So the token function $f_1$ is valid. By induction hypothesis, 
Painter can choose a subset $X_1$ in this round.

In the second phase,   if $G_2$ is a triangle, then $f(c_n)=(4,0)$ since $c_1$ and $c_n$ has a common neighbour $c_j \ne c_2$ on the boundary of $G$.
 Painter colours $c_n$  if $M(c_n) > 0$ and $c_j\notin X_1$. 
 As $c_n$ has four tokens and may lose  two tokens, it can be fully coloured.
 Assume $G_2$ is not a triangle. Painter plays the\dpg on $G_2$ with special vertices $c_1,c_j,c_n$, where  Lister's move $M_2$ is the restriction of $M$ to $V(G_2)$ except that  
 $$M_2(c_j)=\begin{cases}
 1&\text{if}~c_j\in X_1,\\
 0&\text{otherwise}.
 \end{cases}$$
 As $c_n$ is not adjacent to $c_j$   and $c_n$ is not marked in the same round as $c_1$,  Lister's move for the game on $G_2$ is legal.  The token function $f_2$ is obviously valid. Painter has a winning strategy by induction hypothesis and can choose a subset $X_2$.

For the common vertices $c_1,c_j$ of both subgraphs, $c_1$ has only two tokens and  is coloured when it is marked, $c_j$ is marked in the game on $G_2$ if and only if it is coloured the the game on $G_1$. Thus, Painter's move in both games are consistent.  The vertex $c_1$ can only receive defect from $c_2$ because it
plays the role of $c_1$ in both games 
and $c_j$ is not marked in the same round as $c_1$. The vertex $c_j$ plays the role of  $c_2$ in the game on $G_2$, and hence  receives no defect from $V(G_2)\setminus V(G_1)$.  Therefore, Painter's move for $G$ is legal.

\end{proof}

 In the remainder of the proof, we assume that
     $C$ has no chord. This implies that $N(c_1)\cap N(c_n)\cap (V(C)-\{c_2\})=\emptyset$, and hence $f(c_n)=(0,2)$.

 First we choose a subgraph $T$ of $G$, called the the trestle of $G$, as follows.
 A {\em fan} $(P[u,v];q)$ is obtained by joining vertex $q$ to every vertex of path $P[u,v]$ from $u$ to $v$ (it is allowed that $P$ is a single vertex, i.e., $u=v$).
  A {\em trestle} $T=(P;p_1,\cdots, p_{m+1};Q;q_0,\cdots,q_{m})$ is a graph such that
  \begin{itemize}
  \item   $P$ and $Q$ are vertex disjoint paths, where $q_0, q_1,   \ldots, q_m$ are
  distinct vertices of $Q$ occurring in this order (not necessarily consecutive),
  $p_1, p_2, \ldots, p_{m+1}$ are
   vertices of $P$ occurring in this order (not necessarily consecutive), possibly with $p_m=p_{m+1}$, but other $p_j$'s are distinct.
  \item For $i=1,2,\ldots, m$, $(Q[q_{i-1}, q_i];p_i)$ and $(P[p_i, p_{i+1}];q_{i})$ are fans.
  \item Each edge of $T$ is an edge of one of these fans.
  \end{itemize}
 The possible configuration of $T$ is depicted in Figures 4.
 
\begin{figure}[!ht]
	\begin{center}
		\includegraphics[scale=0.25]{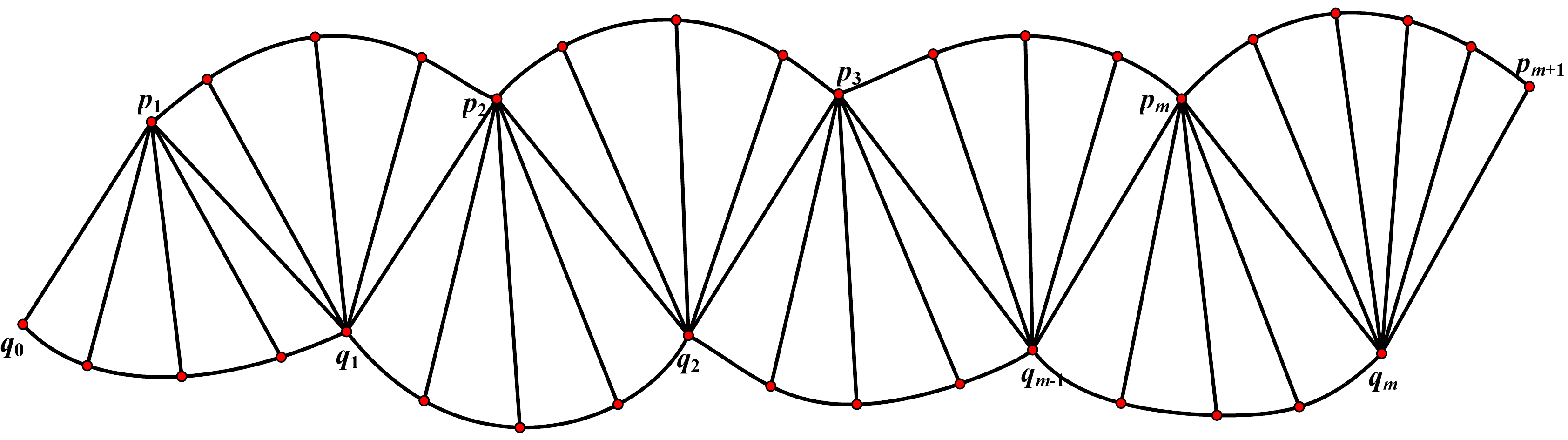}
		\caption{The construction of Trestle.}
		
	\end{center}
\end{figure}

The {\em trestle $T$ of $G$} is the maximal trestle  $T=(P;p_1,\cdots, p_{m+1};Q;q_0,\cdots,q_{m})$ contained in $G$ for which the following hold:
\begin{itemize}
\item $P[p_1, p_{m+1}]$ is a subpath of $C-\{c_1,c_2\}$ starting from $p_1=c_3$, and $Q[q_0, q_m]$ is a path on the boundary of $G-P[p_1,p_{m+1}]$
starting from $q_0=c_2$.
\item Interior vertices of path $Q$ are  not adjacent to any vertex of $C-P$.
\end{itemize}

Let $G_1$ be the subgraph induced by $\{p_i:i=1,2\cdots,m+1\}$. 

If $p_{m+1}=c_n$, then let $G_2=G-P$. 
Otherwise let   $i$ be the maximum index such that $c_i$ is adjacent to $q_m$, and 
 let  $G_2={\rm int}[Q\cup\{q_mc_i\}\cup C[c_{i},c_2]]$.

Let $G_3$ be the subgraph induced by $\bigcup_{i=1}^{m}P(p_i,p_{i+1})$
 and $G_4={\rm int}[C[p_{m+1},c_i]\cup \{c_iq_m,q_mp_{m+1}\}]$.
 If $p_{m+1} \ne c_n$ (as shown in Figure 5(a)),  then  $G_4$ is not empty and $N(q_m)\cap N(p_{m+1})\cap C(p_{m+1},c_i)=\emptyset$. So we just need to ensure that $f_4(p_{m+1})\succeq (0,2)$ for the game on $G_4$.
 If $p_{m+1}=c_n$, then $G_4$ is empty (as shown in Figure 5(b)).

\begin{figure}[ht]
	\begin{center}
		\includegraphics[scale=0.4]{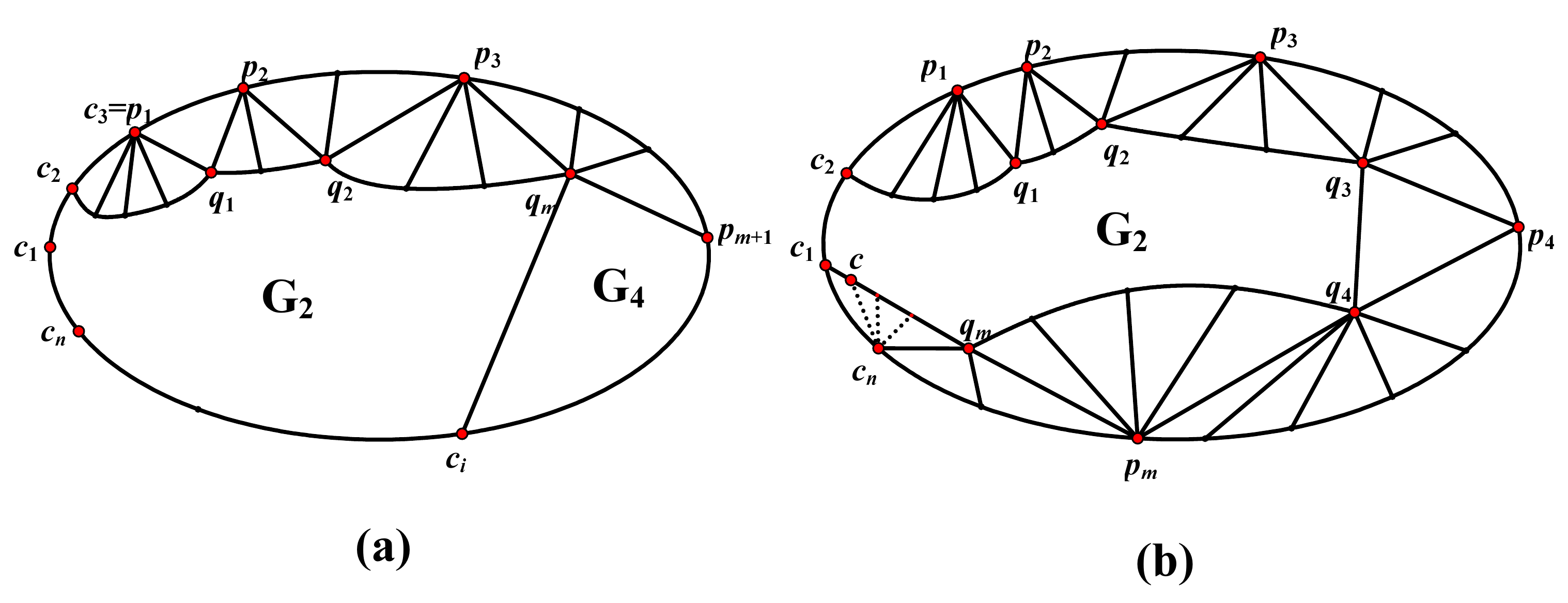}
		\caption{The possible Trestle in $G$.}
		
	\end{center}
\end{figure}

 Painter colours $G_1$ as follows:
 For each $i=1,2,\cdots,m$, let $t_i$ be the number of interior vertices of the path  $P[p_i,p_{i+1}]$
 and if $t_i > 0$, then let  $P(p_i,p_{i+1})=x^i_1x^i_2\cdots x^i_{t_i}$.
 If $t_i>0$, $M(p_i)=2$ and $M(x^i_1)=1$, then we say $p_i$ and $x^i_1$ {\em conflict each other} in this round.
 We consider vertices  $p_1,p_2,\cdots, p_{m+1}$ in order. Note that if $p_{m}=p_{m+1}$, then we consider vertices $p_1,p_2,\cdots, p_m$ in order.  If $p_{m+1} \ne c_n$, then
 colour $p_i$ if and only if all the following hold:
\begin{enumerate}
	\item $M(p_i)=2$.
	\item If $i=1$, then $M(c_2)=0$;  if $i \ge 2$, then $t_{i-1}>0$ or $p_{i-1}$ is not coloured.
	\item This is not the first round in which  $p_i$ and $x^i_1$ conflict each other. 
\end{enumerate}
If $p_{m+1}=c_n$, then  the vertices of $V(G_1)-\{p_{m+1}\}$ are coloured in the same way as above, and $p_{m+1}$ is coloured if and only if $M(p_{m+1})=2$.

For $p_i\in V(G_1)-\{p_{m+1}\}$, $p_i$ may lose two one-dollar token (by (1)),  two tokens due to $p_{i-1}$ or $c_2$ (by (2))  and  one token due to $x^i_1$ (by (3)), and is still a  $\succeq(0,2)$-vertex, so it can be fully coloured. If $p_{m+1} \ne c_n$, then for the same reason, it can be fully coloured. If $p_{m+1}=c_n$, then it is coloured whenever it is marked with a two-dollar token.
As $c_n$ is a $(0,2)$-vertex, it is fully coloured. 

 Moreover,  $p_i$ receives no defect from the vertex in $G_1$, except that if $p_{m+1}=c_n$ and $p_mp_{m+1}\in E(G)$ (or in case that $p_{m-1}p_{m+1}\in E(G)$ if $p_m=p_{m+1}$), then  $p_{m+1}$ and $p_{m}$ (or  $p_{m+1}$ and $p_{m-1}$, respectively) may contribute one defect to each other.  

For the configurations in Figure 5(b), either the other neighbour of $c_1$ on the boundary of $G_2$ is not adjacent to $c_2$ or
the other neighbour of $c_2$ on the boundary of $G_2$ is not adjacent to $c_1$. By symmetry, we may assume 
the other neighbour $c$ of $c_1$ on the boundary of $G_2$ is not adjacent to $c_2$ (as shown in Figure 5(b)). 

In the second phase, Painter plays the\dpg on $G_2$ with   special vertices $c_1,c_2,c_n$ (for the configuration in Figure 5(a)) or with special vertices $c_1,c_2,c$ 
 (for the configuration in Figure 5(b)).
Lister's move $M_2$ is the restriction of
$M$ to $V(G_2)$ except that 
\begin{itemize}
	\item For each $i=1,2,\ldots, m$, for $v \in   N(p_i)-\{q_{i-1},q_i,c\}$, $M_2(v)=0$ whenever $p_i \in X_1$.
	\item For $q_i\ne c$
	$$M_2(q_i)=\begin{cases}
	0&\text{if}~p_{i+1}\in X_1,\\
	\text{max}\{0, M(q_i)-1\}&\text{if}~p_{i+1}\notin X_1 ~\text{and} ~p_i\in X_1.
	\end{cases} 
	$$
	\item For the  configuration  in Figure 5(b), 
	\subitem  if $c\ne q_m$, then for $v\in N(c_n)-\{c_1,q_m,c\} $,     $M_2(v)=0$ whenever $c_n \in X_1$,
	and $M_2(c)=0$ whenever $M(c_1)+M(c_n) \ne 0$;
	\subitem if $c=q_m$, then
	$$M_2(c)=\begin{cases}
	0&~\text{if}~~M(c_1)+M(c_n)\ne 0 ~\\
	\text{max}\{0,M(c)-1\}&\text{ if } ~M(c_1)+M(c_n)= 0 ~\text{ and }~p_m\in X_1.
	\end{cases}$$
\end{itemize}   

For the configuration in Figure 5(a), by the definition of $M_2$, the vertex $q_i$ may lose two tokens and devalue two tokens, and hence $f_2(q_i)\succeq(2,5)$. Each vertex $v\in Q(q_{i-1},q_i)$ may lose two tokens, and  hence $f_2(v)\succeq(0,7)$. By the construction of the trestle,  
$N(c_1)\cap N(c_n)\cap (V(C_2)-\{c_2\})=\emptyset$. So $f_2(c_n)=(0,2)$ is valid at $c_n$.
 Other vertices of $G_2$ obviously have enough tokens, so the token function $f_2$ is valid. 

For the configuration in Figure 5(b), for $i=1,2,\ldots,m-1$, the vertex $q_i$ may lose two tokens and devalue two tokens, and hence $f_2(q_i)\succeq(2,5)$. If $q_m\ne c$, then similarly $f_2(q_m)\succeq(2,5)$ or $f_2(q_m)\succeq(0,7)$ if $p_m=p_{m+1}$. The vertex $c$ is a $(0,9)$-vertex and may lose four tokens, so $f_2(c)\succeq(0,5)$. If $q_m=c$,  then $c$ may lose four tokens and devalue two tokens, and hence $f_2(c)\succeq(2,3)$.  So  $f_2$ is valid at $c$ in both cases. For $i=1,2,\cdots,m$,  the vertex $v\in (N(p_i)-\{q_{i-1},q_i,c\})\cup (N(c_n)-\{c_1,q_m,c\})$ may lose two tokens, so $f_2(v)\succeq(0,7)$. Other vertices of $G_2$ obviously have enough tokens. So the token function $f_2$ is valid.

Thus, Painter can choose a subset $X_2$ by induction hypothesis in both configurations in Figure 5.

In the third phase, for each $i=1,2,\cdots,m$, if $t_i > 0$, then Painter considers the vertices $x^i_{t_i},x^i_{t_i-1},\cdots, x^i_1$ in order.  
For convenience, let $x^i_{t_i+1}=p_{i+1}$. 
For $j=\{1,2,\cdots,t_i\}$, Painter colours $x^i_j$  if and only if all the following hold:
\begin{enumerate}
	\item $q_i\notin X_2$.
	\item  $x^i_{j+1}$ is not coloured.
	\item if $j=1$, then this is not the second round in which  $p_i$ and $x^i_1$ conflict each other. 
\end{enumerate}
 For $j \ge 2$,   $x^i_j$ is a $(2,5)$-vertex and  may lose four tokens due to $q_i$ and $x^i_{j+1}$. So it can be fully coloured.
 Similarly, $x^i_1$ may lose four tokens due to $q_i$ and $x^i_{2}$,  and one $1$-dollar token due to $p_i$,  and is still a $\succeq(1,1)$-vertex.
 So it also can be fully coloured.
 
Observe that the vertex $x^i_j$ receives no defect from $G_3$.

In the fourth phase, Painter plays the\dpg on $G_4$ with special vertices $q_m,c_i,p_{m+1}$, where Lister's move   $M_4$ is the restriction of $M$ to $V(G_4)$, except that   for $v\in \{c_i,q_m\}$,
$$M_4(v)=\begin{cases}
1&\text{if}~v\in X_2,\\
0&\text{otherwise},
\end{cases} $$
and 
$$M_4(p_{m+1})=\begin{cases}
2&\text{if}~p_{m+1}\in X_1,\\
0&\text{otherwise}.
\end{cases} $$ 
Note that $p_{m+1}\ne c_n$, otherwise $G_4$ is empty.  
By the winning strategy in the game on   $G_1$, $p_{m+1}$  is coloured in the round in which $M(p_{m+1})=2$.  So $M_4(p_{m+1}) \le M(p_{m+1})$.
By the definition of $M_4$,  $f_4(p_{m+1})\succeq(0,2)$. The other vertices of $G_4$ obviously have enough tokens.   Thus,  $f_4$ is valid. Painter can choose a a subset $X_4$ by the winning strategy.

The common vertices  $q_m,c_i$ of $G_4$ and $G_2$ (respectively, the common vertex $p_{m+1}$ of $G_4$ and $G_1$)
are coloured in the game on $G_4$ if and only if they are coloured in the game on $G_2$ (respectively, in $G_1$).  
So Painter's move for all the subgraphs are consistent.  

Moreover, as $q_m,c_i$ play the roles of $c_1,c_2$ in the game on $G_4$,  they   receive  no defect from $G_4 -\{q_m,c_i\}$.  So $q_m,c_i$ satisfy Rule (P1). If $p_{m+1}\ne c_n$, then   $p_{m+1}$  is the common vertex of $G_4$ and $G_1$, and  receives no defect from $G_i$ for $i=1,2,3$.  If $p_{m+1}=c_n$, then $p_{m+1}$ receives at most one defect from $p_m$. Moreover, if $p_m=p_{m+1}$, then by the same reason, $p_{m+1}$  receives at most  one defect 
    from $ p_{m-1}$. So $p_{m+1}$ satisfies Rule (P1).
    
For $i=1,2,3,\cdots, m$, the vertex $p_i$ is coloured in the round in which $M(p_i)=2$ and  receives at most  one defect from $\{ q_i,x^i_1\}$ (note that
$q_i$ and $x^i_1$ are not coloured in the same round). In case, 
  $p_{m+1}=c_n$  and $p_{m}p_{m+1}\in E(G)$, then $p_m$   may receive one defect from $p_{m+1}$. As $q_m$  and $p_{m+1}$ are not coloured in the same round, 
$p_m$  receives at most  one defect from $\{ q_m,p_{m+1}\}$.  So vertices of $\{p_i: i=1,2, \cdots, m\}-\{p_{m+1}\}$ satisfy Rule (P1).

   Note that if $p_i$ and $q_i$ are coloured in the same round, then $q_i$ is marked with a one-dollar token in the game on $G_2$, 
   and hence receives no defect from $G_2$.  If $p_i$ and $x^i_{1}$ are coloured in the same round, then $p_i$ and $x^i_1$ do not conflict
   each other in this round. So $x^i_1$ is coloured with a $2$-dollar token in this round. Moreover, $x^i_1$   receives no defect from $G_3$. So $x^i_1$ satisfies Rule (P1). Therefore,    Painter's move is legal for the game on $G$.

\section*{References}

\bibliographystyle{siam}

\bibliography{ohba}

\begin{thebibliography}{10}
\expandafter\ifx\csname url\endcsname\relax
  \def\url#1{\texttt{#1}}\fi
\expandafter\ifx\csname urlprefix\endcsname\relax\def\urlprefix{URL }\fi

\bibitem{ATV97}
N.~Alon, Z.~Tuza, M.~Voigt, Choosability and fractional chromatic numbers,
  Discrete Math. 165/166 (1997) 31--38, graphs and combinatorics (Marseille,
  1995).

\bibitem{CZ}
T.~P. Chang, X.~Zhu, On-line $3$-choosable planar graphs, Taiwanese Journal of
  Mathematics 16 (2012) 511--519.

\bibitem{CCW1986}
L.~J. Cowen, R.~H. Cowen, D.~R. Woodall, Defective colorings of graphs in
  surfaces: partitions into subgraphs of bounded valency, J. Graph Theory 10
  (1986) 187--195.

\bibitem{Kierstead2010}
W.~Cushing, H.~A. Kierstead, Planar graphs are $1$-relaxed, $4$-choosable,
  European J. Combin. 31 (2010) 1385--1397.

\bibitem{Eaton1999}
N.~Eaton, T.~Hull, Defective list colorings of planar graphs, Bull. Inst.
  Combin. Appl 25 (1999) 79--87.

\bibitem{ERT79}
P.~Erd{\H{o}}s, A.~L. Rubin, H.~Taylor, Choosability in graphs, in: Proceedings
  of the {W}est {C}oast {C}onference on {C}ombinatorics, {G}raph {T}heory and
  {C}omputing ({H}umboldt {S}tate {U}niv., {A}rcata, {C}alif., 1979), Congress.
  Numer., XXVI, Utilitas Math., Winnipeg, Man., 1980.

\bibitem{Gutowski2011}
G.~Gutowski, Mr. paint and mrs. corrector go fractional, Electron. J. Combin.
  18~(1) (2011) Research Paper 140.

\bibitem{GGHZ2016}
G.~Gutowski, M.~Han, T.~Krawczyk, X.~Zhu, Every planar graph is $3$-defective
  $3$-paintable, manuscript.

\bibitem{HZ2014}
M.~Han, X.~Zhu, Locally planar graphs are $5$-paintable, Discrete Math. 338
  (2015) 1740--1749.

\bibitem{HZ2015}
M.~Han, X.~Zhu, Locally planar graphs are $2$-defective $4$-paintable, European
  Journal of Combinatorics 54 (2016) 35--50.

\bibitem{HRS1973}
A.~Hilton, R.~Rado, S.~S.H., A $(<5)$-colour theorem for planar graphs, Bull.
  London Math. Soc. 5 (1973) 302--306.

\bibitem{Schauz2009}
U.~Schauz, Mr. {P}aint and {M}rs. {C}orrect, Electron. J. Combin. 16~(1) (2009)
  Research Paper 77, 18.

\bibitem{SKI1999}
R.~{\v S}krekovski, List improper colourings of planar graphs, Combin.
  Probability and Computing 8 (1999) 293--299.

\bibitem{Tho1994}
C.~Thomassen, Every planar graph is {$5$}-choosable, J. Combin. Theory Ser. B
  62~(1) (1994) 180--181.
\newline\urlprefix\url{http://dx.doi.org/10.1006/jctb.1994.1062}

\bibitem{Zhu2016}
X.~Zhu, Multiple list colouring of planar graphs, manuscript.

\end{thebibliography}

\end{document}